\documentclass[review,onefignum,onetabnum]{siamart171218}

\usepackage{lipsum}
\usepackage{amsfonts}
\usepackage{graphicx}
\usepackage{epstopdf}
\usepackage{algorithmic}
\ifpdf
  \DeclareGraphicsExtensions{.eps,.pdf,.png,.jpg}
\else
  \DeclareGraphicsExtensions{.eps}
\fi

\newsiamremark{remark}{Remark}
\newsiamremark{hypothesis}{Hypothesis}
\crefname{hypothesis}{Hypothesis}{Hypotheses}
\newsiamthm{claim}{Claim}

\headers{Read Between the Hyperplanes: Extensions to Randomized Kaczmarz}{J. Nguyen, O. Presnyakov, and A. Radhakrishnan}

\title{Read Between the Hyperplanes: \\On Spectral Projection and Sampling \\Approaches to Randomized Kaczmarz\thanks{Advised by Dr. Anna Ma (University of California, Irvine)}}

\author{James A. Nguyen\thanks{Department of Mathematics, University of California, Irvine.}
\and Oleg Presnyakov\footnotemark[2]
\and Adityakrishnan Radhakrishnan\footnotemark[2]
}

\usepackage{amsopn}

\makeatletter
\newcommand*{\addFileDependency}[1]{
  \typeout{(#1)}
  \@addtofilelist{#1}
  \IfFileExists{#1}{}{\typeout{No file #1.}}
}
\makeatother

\newcommand*{\myexternaldocument}[1]{%
    \externaldocument{#1}%
    \addFileDependency{#1.tex}%
    \addFileDependency{#1.aux}%
}

\newenvironment{mydef}[1][Definition]{%
  \refstepcounter{equation}
  \par\noindent\textbf{#1~\theequation.}\ }{\par}

\ifpdf
\hypersetup{
  pdftitle={Read Between the Hyperplanes: On Spectral Projection and Sampling Approaches to Randomized Kaczmarz},
  pdfauthor={J. Nguyen, O. Presnyakov, and A. Radhakhrishnan}
}
\fi

\myexternaldocument{ex_supplement}

\begin{document}
\nolinenumbers

\maketitle

\begin{abstract}
  Among recent developments centered around Randomized Kaczmarz (RK), a row-sampling iterative projection method for large-scale linear systems, several adaptions to the method have inspired faster convergence. Focusing solely on ill-conditioned and overdetermined linear systems, we highlight inter-row relationships that can be leveraged to guide directionally aware projections. In particular, we find that improved convergence rates can be made by (i) projecting onto pairwise row differences, (ii) sampling from partitioned clusters of nearly orthogonal rows, or (iii) more frequently sampling spectrally-diverse rows.
\end{abstract}



\section{Introduction}

The size of computational models and applications continues to grow at an exponential rate. Thus, we seek faster, more robust algorithms to handle the large-scale linear systems that depend on them. In the past decade, two prevalent examples, Randomized Kaczmarz (RK) and Sampling Kaczmarz-Motzkin (SKM), have seen drastic improvement in contributions that utilize an added element of randomization. The aforementioned methods add an element of randomization and partial greediness to promote faster convergence \cite{haddock2020greedworksimprovedanalysis}. These methods are row-iterative projection methods for solving the aforementioned large-scale linear systems $Ax = b$, where $x^*$ is the intended solution in the set of feasible solutions $x = \{x \in \mathbb R^m : Ax \leq b\}$. RK randomizes the sampling of the existing Kaczmarz projection method; if $A \in \mathbb R^{m\times n}$, the method updates the current iterate as follows:
\begin{align}
    x_{k+1} = x_k - \frac{\langle a_i, x_k \rangle - b_i}{||a_i||^2} a_i,
\end{align}
\noindent where $a_i$ is the $i^{th}$ row, randomly sampled from $1, \dots, m$. Traditionally, it is common to use a sampling distribution that selects rows with probability proportional to its squared Euclidean norm $||a_i||^2$; this extension allows us to guarantee, within reasonable probability, convergence of the $k^{th}$ iterate $x_k$ to $x^*$ that is at least linear in expectation: 
\begin{align}\label{expectedConvergence}
    \mathbb E(||\varepsilon_k||^2) \leq \left(1- \frac{\sigma_{\min}^2(A)}{||A||_F^2} \right)^k||\varepsilon_0||^2,
\end{align}

\noindent where $\epsilon_k$ is the $k^{th}$ error vector $x_k-x^*$ \cite{strohmer2009randomized}. Similarly, the SKM algorithm takes a partially greedy approach to sampling, promoting faster iterative convergence at a higher computational cost. Given a greediness parameter $\beta \in \mathbb N$, this algorithm samples without replacement a set $\tau_k$ with $|\tau_k| = \beta$ from the system $A \in \mathbb R^{m\times n}$, then it projects onto the row $a_i \in \tau_k$ that maximizes the residual $\langle a_i, x_k\rangle -b_i.$ By greedily choosing a row from $\tau_k$ at each iterate $k$, the convergence in expectation improves drastically \cite{deloera2019samplingkaczmarzmotzkinalgorithmlinear}. \\

Note that the bound for the $k^{th}$ expected error of either algorithm is dependent on the scaled conditioning of the matrix. Consequently, we would expect that the algorithms would have weaker theoretical guarantees for ill-conditioned linear systems. Thus, we seek to explore several methods for sampling and clustering to address RK and SKM's limitations under these conditions. For the remainder of this study, we restrict our attention to the ill-conditioned case for overdetermined consistent linear systems.

Throughout this paper, we suggest three different methods to improve the convergence rate and approximation error of RK and SKM: (i) a transformed linear feasibility problem that considers pairwise row differences in Section \ref{sec:pairwise}, (ii) the construction of a coreset to use partitions of the system to compute its solution in Section \ref{sec:clustering}, and (iii) an increased sampling frequency of rows corresponding to most underrepresented singular values of the system in Section \ref{sec:singular_directions}.

\section{Pairwise Comparisons for Binary Classification}
\label{sec:pairwise}

In the following section, we offer an alternative solution to solving overdetermined binary classification problems using row-iterative methods. We propose a procedure to re-imagine a binary classification problem as a linear feasibility one, and we aim to leverage variations of RK to iteratively converge toward a solution. To transform the problem, we define the following Hadamard product operation to scale each row of the system by the sign of its label.
\begin{mydef}
    [Hadamard Product]
    Let $A, B \in \mathbb R^{m\times n}$. Then, we take the Hadamard product to be $\odot: \mathbb R^{m\times n} \times \mathbb R^{m\times n} \to \mathbb R^{m\times n}$ defined by:
    \[
    (A\odot B)_{ij} = A_{ij}B_{ij},
    \]
    the element-wise matrix product.
\end{mydef}

\noindent We consider a consistent system of inequalities $Ax^* \leq b$ for some matrix $A$ and an unknown intended solution vector $x^*$ such that
\[
b_i = \begin{cases}
    -1, & (Ax^*)_i < 0 \\
    1, & (Ax^*)_i \geq 0
\end{cases}.
\]
To setup the linear feasibility problem, we obtain the matrix $B \in \mathbb R^{m\times n}$, where $B_{:j} = b$ for all $j = 1, \dots, n$. Using this matrix, along with the Hadamard product, we define:
\[
A' := -B \odot A
\]
and the zero vector $b' = 0$. Thus, we end up with the linear feasibility problem $A'x^* \leq b'$.\\

\noindent Recall that the original goal of the problem is to find a vector $x$ such that 
\[
\text{sign}(\langle A_i, x\rangle) = \text{sign}(\langle A-i, x^* \rangle) = b_i.
\]
Note that the following conditions are equivalent:
\begin{align}
    b_i &\times \text{sign}(\langle A_i,x\rangle) > 0  \label{problem:binary_classification}\\
    &\iff (B_iA_i)x > 0 \nonumber\\
    &\iff ((-B \odot A)x)_i < 0 \nonumber\\
    & \iff (A'x)_i < b'_i \nonumber
\end{align}
Thus, we can solve a binary classification problem for an overdetermined linear system by using row-iterative methods such as RK and SKM, as long as we transform it by the described procedure. 

\subsection{Approach} \label{subsection:pairwise_approach}
Although it would be valid to simply apply the RK or SKM algorithm to the resulting linear feasibility problem, we suggest that additional informative comparisons can be made between rows in the system. To test this hypothesis, we consider extracting from a matrix $A'$, the pairwise differences matrix $P \in \mathbb R^{\binom {m}{2} \times n}$ defined by:
\[
P_{h,:} = A'_i - A'_j, 
\]
where $h \in \{1, \dots, \binom m2$\} is a unique combination of indices $i \neq j$, and $i,j \in \{1, \dots, m\}$. If we row bind $A'$ and $P$, we obtain a combined matrix $A'\cup P \in \mathbb R^{m + \binom m2 \times n}$ with additional comparisons.  Since the rows of $A'$ have been scaled by the signs of $Ax^*$, we have an similar problem to \ref{problem:binary_classification}. We can again proceed to using RK or SKM to find a vector $x \in \mathbb R^{m + \binom m2}$ that satisfies $(A'\cup P)x < 0$.\\

\noindent If we define $y = x_{1:m,:}$, we hypothesize that the $k^{th}$ iterate $y^{(k)}$ converges to $x^*$. Furthermore, we hypothesize that it will converge faster than the iterate $x'^{(k)}$, obtained by applying RK or SKM to find $x' \in \mathbb R^m$ satisfying $A'x' < 0$. \\

\subsection{Methods} 

To analyze this relationship, we run numerical experiments on 100 random matrices, obtained by taking reversing the Singular Value Decomposition of two random orthonormal matrices $U, V$ and a diagonal matrix $S$ an exponentially scaled singular values to replicate an ill-conditioned system. \\

\noindent For the following experiments, we employ the SKM algorithm with a greediness parameter $\beta = 3$. 

\begin{algorithm}
\caption{SKM for linear feasibility}\label{alg:skm}
\begin{algorithmic}
\STATE \textbf{procedure} SKM($A$, $b$, $x_0$, $x^*$, $\beta$, $\lambda$, $K$)
\STATE \hspace{0.1in} $k=1$
\STATE \hspace{0.1in} \textbf{while $k < K$}
\STATE \hspace{0.2in} Choose a sample of $\beta$ rows, $\tau_k$, from $A.$
\STATE \hspace{0.2in} $\displaystyle t_k := \operatorname*{argmax}_{i \in \tau_k} a_i^Tx_{k-1} - b_i$
\STATE \hspace{0.2in} $x_k = x_{k-1} - \lambda(a_{t_k}^Tx_{k-1} - b_{t_k})^+a_{t_k}$
\STATE \hspace{0.2in} $k = k+1$
\STATE \hspace{0.1in} \textbf{return} $x_k$
\STATE \textbf{end procedure}
\end{algorithmic}
\end{algorithm}

\noindent We apply the algorithm to the matrix $(A' \cup P)$ described in subsection \ref{subsection:pairwise_approach}, and we compare the average approximation error, Chebyshev error, and binary classification accuracy of the algorithm performed using three different sampling schemes: 
\begin{enumerate}
    \item first $m$ rows sampled uniformly, remaining $\binom m2$ rows not sampled, 
    \item all $m + \binom m 2$ rows weighted uniformly,
    \item first $m$ rows not sampled, remaining $\binom m2$ rows sampled uniformly.
\end{enumerate} 

\noindent We define the Chebyshev error to be the distance from the iterate vector to the Chebyshev center of the feasible region.
\begin{mydef}
    [Chebyshev Center of a Feasible Region] The point within a feasible region that maximizes the volume of a feasible sphere centered around it.
\end{mydef}

\subsection{Results}

\begin{figure}[h]
    \centering
    \includegraphics[width=0.32\linewidth]{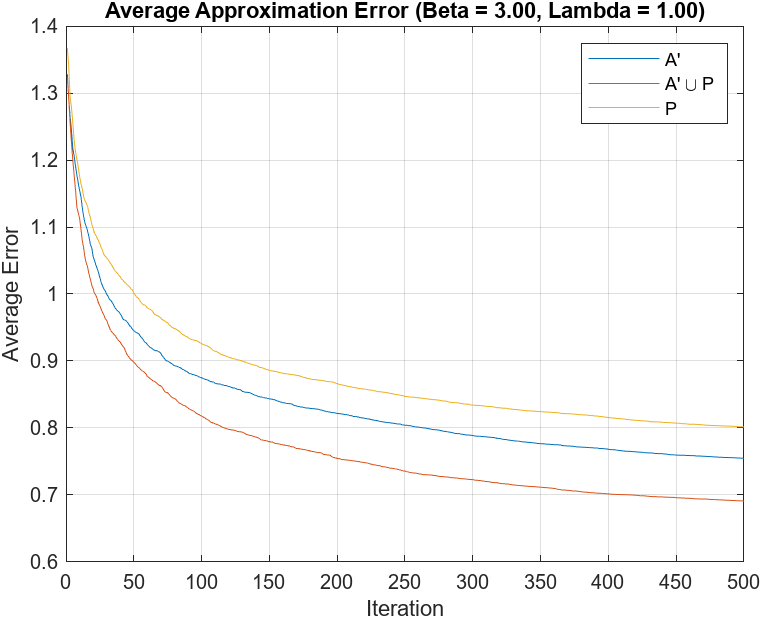}
    \includegraphics[width=0.32\linewidth]{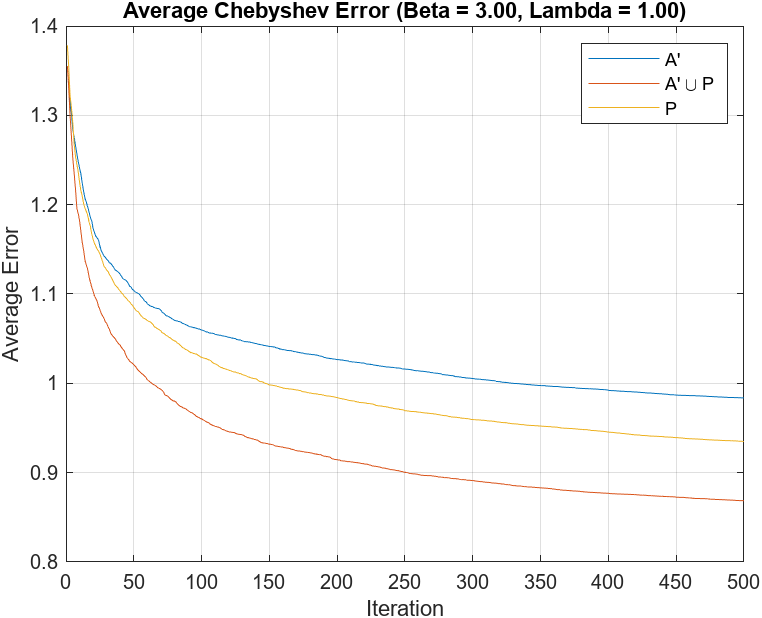}
    \includegraphics[width=0.32\linewidth]{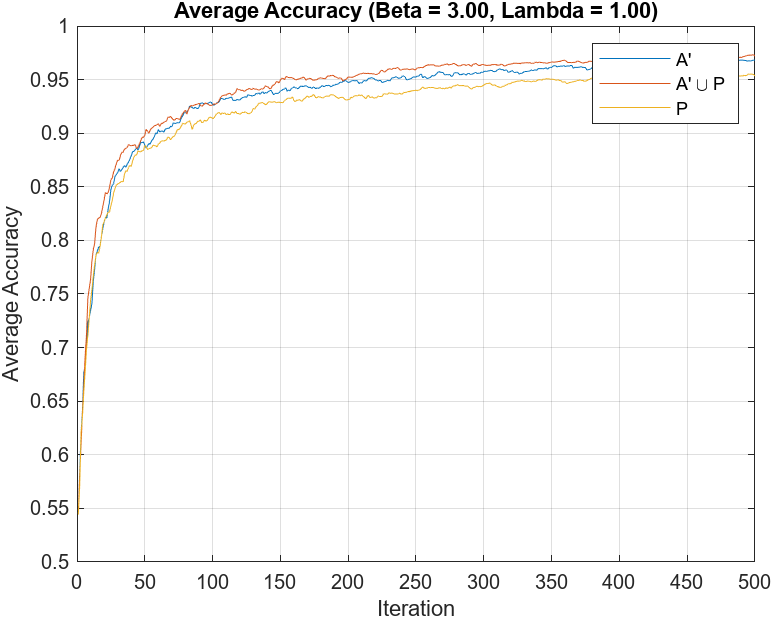}
    
    \caption{Quantitative results for $A' \in \mathbb R ^{240\times 12}$, averaged over 100 trials: (i) Iteration vs. Average Approximation Error, (ii) Iteration vs. Average Chebyshev Error, and (iii) Iteration vs. Average Accuracy (left to right). All results obtained by applying SKM using $\beta = 3$, with no over-projection parameter $(\lambda = 1)$. The above plots feature metrics associated with the sampling schemes 1, 2, and 3, referenced by $A', A'\cup P, $ and $P$, respectively.}
    \label{fig:placeholder}
\end{figure}

The numerical experiments suggest that sampling scheme 2 retains similar success in producing an accurate binary classifier; however, it reveals that additionally sampling pairwise differences promotes a reduced average approximation error and Chebyshev error. Although it seems that sampling a pairwise difference row could be more informative, we also see that strictly sampling these rows produces a less accurate binary classifier with greater approximation error. Thus, it should not be assumed that these comparisons are strictly more informative of the true solution to the system.

\section{Clustering}
\label{sec:clustering}

Another central question we investigated concerns the extraction of a compact coreset from an over-determined linear system. We conjecture that any system of size  $m × n$ admits a coreset of cardinality $O(c*n)$, where  $c$ is a small constant independent of $m$. Identifying such a subset would enable existing algorithms to operate on a substantially reduced system while preserving the original problem’s solution quality. Our approach is to partition the matrix  $A$ into clusters of similar rows, treating each cluster as an approximately redundant group whose representative captures the essential structure of that subset. This clustering strategy allows for parallel processing and lowers the overall computational cost without significant loss of accuracy. Throughout our experiments, we use the Kaczmarz algorithm as the baseline solver for evaluating convergence speed and approximation fidelity. The following section elaborates on the theoretical intuition and presents the algorithms that implement this idea.

\subsection{Existence of a ``Good" Cluster}
To obtain a quick proof–of–concept that a small \emph{inner‑product coreset}
can solve an over‑determined linear system almost as accurately as the full
data set, we ran the following procedure on a synthetic ill-conditioned matrix.

\begin{enumerate}
    \item \textbf{Data generation.}
        For most of our tests we generated \emph{ill‑conditioned matrices}, which usually constitute the worst‑case scenario for iterative algorithms. Consequently, insights gained on such matrices tend to generalize well to arbitrary random ones. To build each ill‑conditioned matrix we employed an SVD‑based procedure: two random orthogonal factors were sampled, and the singular values were manually assigned so that the ratio between the largest and the smallest was on the order of $10^{7}$.
              
    \item \textbf{Row scoring.}  
          For every row \(a_i\) compute the absolute inner product  
          \(\;s_i = \bigl|\langle a_i,\,x^\star\rangle\bigr|\).  
          Large scores indicate rows whose directions are highly aligned with
          \(x^\star\).
    \item \textbf{Coreset extraction.}  
          Sort the scores and keep the \(k= c*n\) rows with the \emph{smallest}
          absolute inner products (we use the heuristic factor \(c=2\)).  The
          retained rows form a sub‑matrix \(B\) with the matching right‑hand
          side \(b_B\).
    \item \textbf{Least‑squares solves.}
          \begin{itemize}
              \item Full system:
                    \(\displaystyle
                         x_A=\arg\min_x \|A\,x-b\|_2.\)
              \item Coreset system:
                    \(\displaystyle
                         x_B=\arg\min_x \|B\,x-b_B\|_2.\)
          \end{itemize}
          Both minimizers are obtained with \texttt{numpy.linalg.lstsq}.
    \item \textbf{Evaluation.}  
          Measure the relative error
          \(\|x_B-x_A\|_2/\|x_A\|_2\) and compare residual norms
          \(\|A x_B-b\|_2\) and \(\|B x_B-b_B\|_2\).
\end{enumerate}

\begin{figure}[h!]
    \centering
    \begin{tabular}{cc}
        \includegraphics[width=0.45\linewidth]{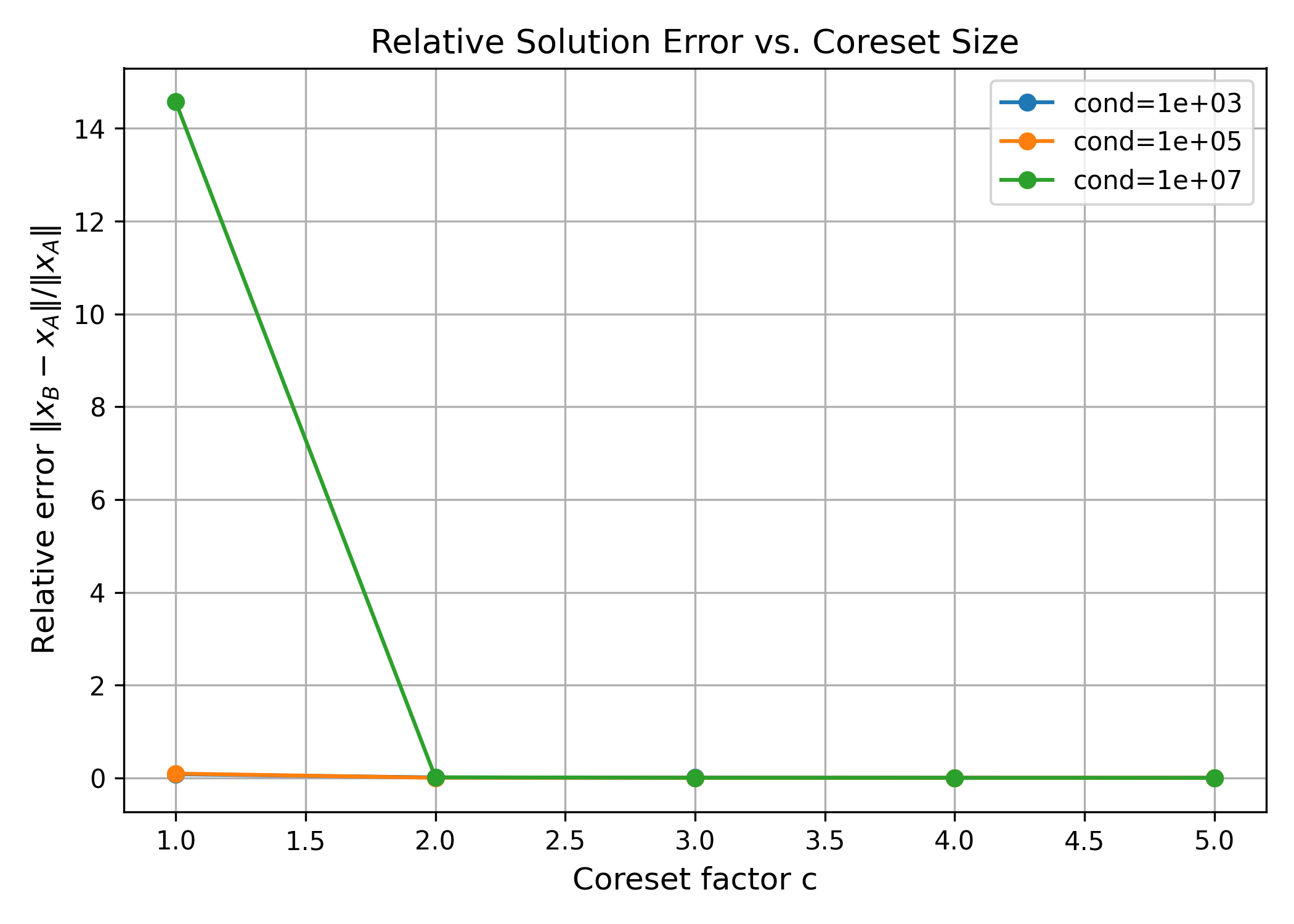} &
        \includegraphics[width=0.45\linewidth]{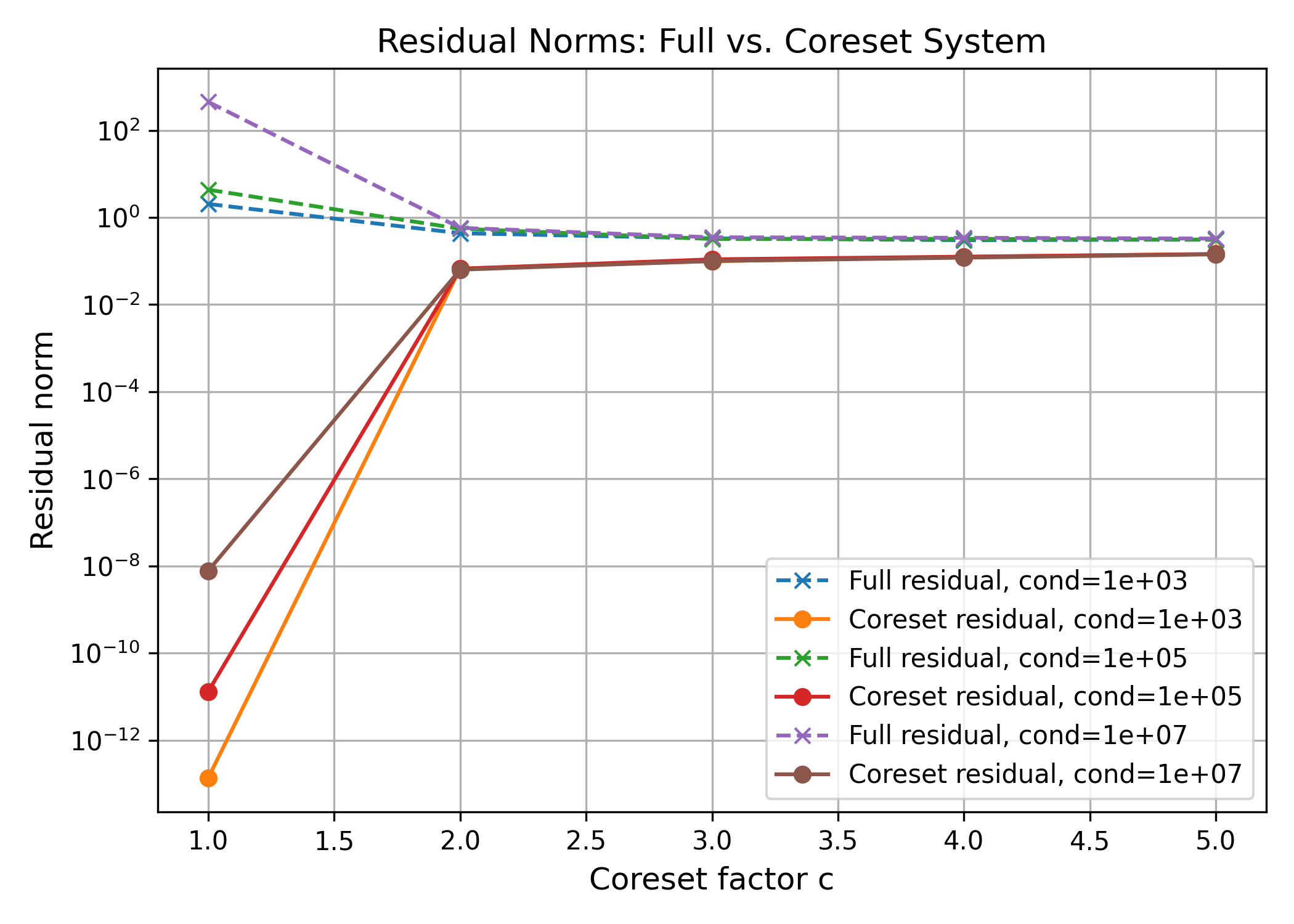}
    \end{tabular}
    \caption{Comparison of coreset performance under varying condition numbers.
    (\textbf{Left}) Relative solution error between the full least-squares solution
    $x_A$ and the coreset solution $x_B$ as a function of the coreset factor $c$.
    Increasing $c$ consistently improves the approximation quality, even for
    ill-conditioned systems ($\kappa(A)=10^7$). 
    (\textbf{Right}) Residual norms for the full and coreset systems plotted on a
    logarithmic scale. Dashed lines denote residuals $\|A x_B - b\|$, while solid
    lines correspond to $\|B x_B - b_B\|$. The stability of residual magnitudes
    across condition numbers indicates that the reduced system preserves the
    dominant geometry of the original problem.}
    \label{fig:coreset_results}
\end{figure}

The results in Figure~\ref{fig:coreset_results} confirm that compact coresets can retain the fundamental subspace geometry of the original system even under extreme ill-conditioning. As the coreset factor $c$ grows, the relative error between $x_B$ and $x_A$ decreases approximately monotonically, suggesting that the selected rows form an increasingly accurate projection of the column space of $A$. The near-constant residual magnitudes across condition numbers indicate that the coreset preserves the dominant left-singular subspace, implying that $B^\top B \approx A^\top A$ in spectral norm up to a small perturbation. These empirical results align with the intuition that low-rank or leverage-based clustering implicitly enforces subspace embedding guarantees, thereby maintaining both conditioning and numerical stability while drastically reducing the system size.

\subsection{Algorithms}
In this section we will talk about all the algorithms we've been trying to utilize in order to achieve the best clustering. At the end of the section we will provide plots with all clustering methods to compare its result on a real data.

\subsubsection{$\epsilon$-cover net Algorithm}
The first algorithm is based on an $\epsilon$-cover net idea. For visualization, imagine each row as a vector in an $n$-dimensional vector space that connects the origin and the point defined by coefficients of a matrix. Then we are trying to understand how each row acts on the other vectors in our space with respect to an inner product. It is not hard to imagine that if points are close to each other, then they act similarly. As a result, we decided to use partition based on how far from each other points are. For details, refer to Algorithm \ref{alg:cluster}.

\begin{algorithm}
\caption{Epsilon-Cover Net Algorithm}
\label{alg:cluster}
    \begin{algorithmic}
        \STATE Assign $a_1$ to Cluster $S_1$.
        \FOR{each $a_k$, $k>1$} 
        \IF{$\exists ~ S_i$ such that $\Big\langle a_k, \tfrac{1}{p}\sum_{j=1}^{p} s_j \Big\rangle < \varepsilon$}
        \STATE Append $a_k$ to $\underset{i}{\arg\min} \; S_i$
        \ELSE
        \STATE Initialize new cluster, append $a_k$.
        \ENDIF
        \ENDFOR
    \end{algorithmic}
\end{algorithm}

\noindent After clustering, we project iterates onto one random row from each cluster. Then pick the ``best'' cluster, $S^*$, where the criterion for ``best'' is being the most orthogonal cluster to the current iterate (which is our best approximation for $x^*$), i.e., we let $S^*$ be cluster $S_i$ where $$ i := \operatorname*{argmin}_{i}|\langle x_k, \frac{1}{p}\sum_{i=1}^{p}s_i\rangle|.$$

The algorithm’s performance hinges on the chosen coverage radius $\epsilon$. If $\epsilon$ is set too small relative to the spacing of the vectors, clusters become sparse and uninformative; if too large, nearly all rows collapse into a single cluster. Hence, our next objective is to identify an $\epsilon$ that yields a well‑balanced matrix partition. In addition, because we track the mean projection in every cluster, we can dynamically reassign rows to different clusters once the current partition stops producing error reductions in further iterations.

\subsubsection{Online cluster updating}
Another idea is essentially taken from a well-known K-means algorithm. Since we know that the Karzmarz algorithm will eventually converge to a solution, we were trying to find an optimization that will be updated during the iterations. Each fixed number of iterations we were trying to reduce the matrix $A$ that we were working with. Particularly, for an $m \times n$ matrix:
\begin{enumerate}
    \item Start with a cluster of the $2n$ best rows from the first active set which is the entire matrix of $m$ rows.
    \item Every $(2^k\times100)$-th iteration ($k=0,1,...$), reduce the size of the active set by a factor of 2 ($\frac{m}{2}, \frac{m}{4}, ...$) and stop at $4n$ rows.
    \item For each active set, we recluster and choose the $2n$ best rows based on the current iterate.
\end{enumerate}

\noindent The ``best'' row is defined as the one most orthogonal to the current iterate:
\[
    i \;:=\;\underset{i \in [m]}{
    \arg\min}
    \left( a_i^{\mathsf T} x_{k-1} \;-\; b_i \right).
\]

\noindent With each iteration, the rows we discard contribute less and less: their projections cease to reveal new information about the solution subspace. Although the exact pattern of elimination depends on the particular right‑hand side, the sequence of discarded sets still carries useful structure. We can harness this history to prospectively form clusters whenever the same matrix must be solved against multiple right‑hand sides.

\subsection{Experimental Evaluation on Gaussian vs.\ Ill-conditioned matrices}
Preliminary experiments on Gaussian matrices did not reveal substantial performance differences between the considered algorithms. In this setting the underlying geometry is ``well-behaved'': random directions are in general position, and each projection step in the Kaczmarz-type methods almost surely yields a sufficiently informative update. As a consequence, all variants exhibit similarly fast convergence on Gaussian data and this regime does not clearly differentiate their behavior.

\begin{figure}[h!]
    \centering
    \includegraphics[width=\linewidth]{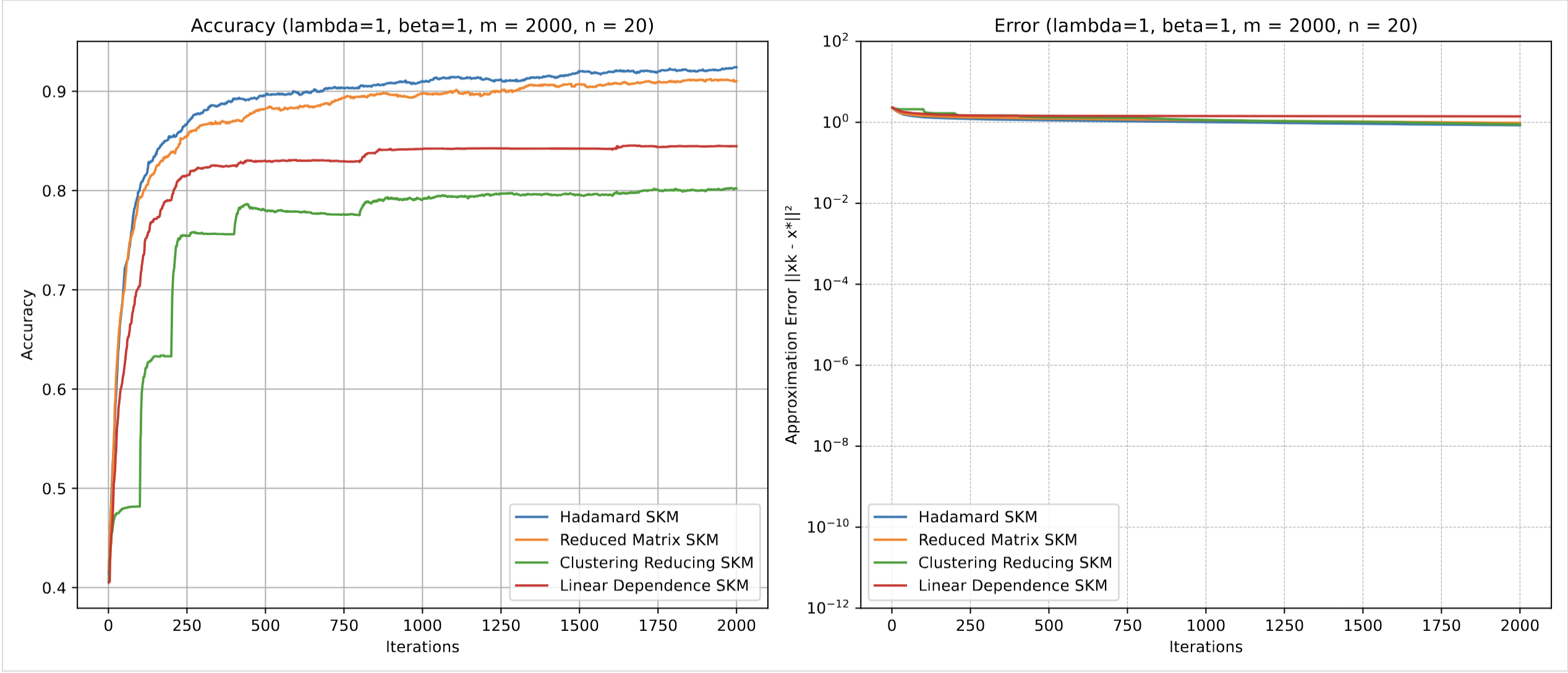}
    \caption{Convergence of SKM variants on an ill-conditioned linear system with $(m,n) = (2000,20)$ and condition number $\kappa(A)\approx 10^{7}$ over $2000$ iterations. The blue curve corresponds to the original Hadamard SKM, the orange curve to the reduced-matrix (proof-of-concept) variant using a coreset of size $c=5$, the green curve to the clustering-based reduction (``$\varepsilon$-cover'') approach, and the red curve to the online clustering / linear-dependence reduction scheme. The figure illustrates that the reduced-matrix variant closely tracks the baseline, while clustering-based variants exhibit slower convergence and larger approximation error.}
    \label{fig:ClusterRun}
\end{figure}

To obtain a more discriminative setting, we therefore evaluated all methods on a synthetically generated ill-conditioned matrix with fixed condition number on the order of $10^{7}$. In Figure~\ref{fig:ClusterRun} we report the convergence of four stochastic Kaczmarz Motzkin (SKM) variants over $2000$ iterations on a system of size $(m,n) = (2000, 20)$. The ``reduced matrix'' (proof-of-concept) variant uses a constant $c = 5$, corresponding to a coreset of about $5\%$ of the rows. In this experiment, the reduced matrix method achieves accuracy and approximation error comparable to the full-matrix Hadamard SKM baseline, demonstrating that a suitably selected subset of rows can solve the linear system nearly as well as the original system. However, none of the modified algorithms matches the performance of the original Hadamard SKM, indicating a nontrivial trade-off between reduction and convergence speed in the ill-conditioned regime.

\begin{table}[h!]
    \centering
    \label{tab:skm_results}
    \begin{tabular}{|lcc|}
        \hline
        \textbf{Method}             & \textbf{Accuracy (\%)} & \textbf{Approx. Error} \\
        \hline
        Hadamard SKM (original)            & 92.42                  & 0.8486 \\
        Reduced Matrix SKM (proof-of-concept) & 91.03               & 0.9697 \\
        Clustering-Reduced SKM             & 80.20                  & 0.8848 \\
        Linear-Dependence SKM              & 84.47                  & 1.4024 \\
        \hline
    \end{tabular}
        \caption{Final performance of SKM variants after $2000$ iterations on an ill-conditioned system with $\kappa(A)\approx 10^{7}$ and parameters $\lambda = 1$, $\beta = 1$. Accuracy is reported in percent; approximation error is measured as the residual norm (lower is better).}
\end{table}

\newpage

\section{Convergence along Singular Directions}
\label{sec:singular_directions}

\subsection{Background}
The latest avenue of our exploration emphasizes a direction-aware approach to RK algorithms. Whilst previous approaches focus on measures of orthogonality, we redirect that focus to analyzing the singular directions present in our system. Our analysis is inspired by recent developments revealing that the error vector converges along the system's smallest singular vector \cite{steinerberger2021}. Steinerberger's analysis of the convergence of the error vector suggests that RK -- under traditional sampling techniques -- struggles to iterate in the directions least represented by the data.  \\

\begin{figure}[h!]
    \centering
    \includegraphics[width=1\linewidth]{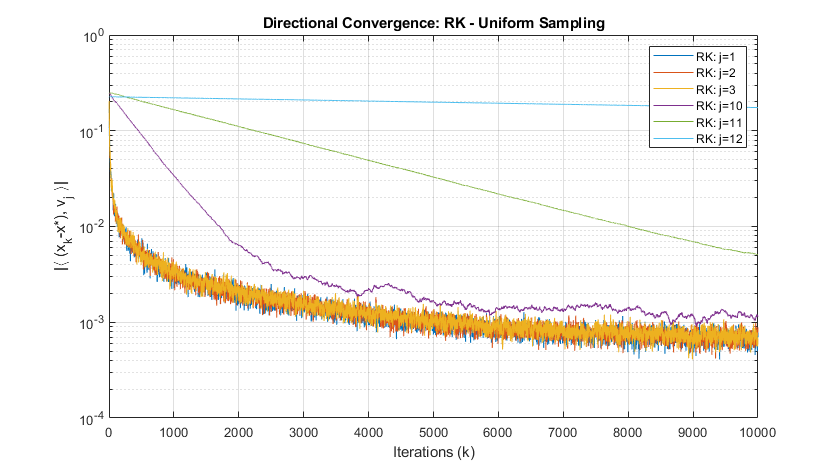}
    \caption{$A \in \mathbb R^{240 \times 12}$ is constructed as described in Section \ref{sub:numericalExperimentation}. The figure demonstrates the convergence of $|\langle x_k-x^*, v_j\rangle |$ for each singular vector $v_j ~(j = 1,\dots, n)$ at each iteration $k$. It demonstrates the ordinal spectrum of convergence by decreasing singular value.} 
    \label{fig:SpectralConvergence}
\end{figure}

\noindent We see that Figure \ref{fig:SpectralConvergence} confirms  Steinerberger's results from \cite{steinerberger2021} and suggest the component of the error vector $\varepsilon_k$ in the direction of the largest singular vectors converge to 0 most quickly.

\subsection{Approach}
To compensate for the under-representation of small singular vectors in the trajectory of our iterate $x_k$, we seek to adapt our sampling distribution to select with higher probability rows that have a greater component in the direction of the least singular vector. Let $a_i$ be the $i$th row of $A \in \mathbb R^{m \times n}$. Then, $a_i$ is some linear combination of the $n$ singular vectors:
\begin{align}
    a_i = \sum_{j=1}^n c^{(i)}_j\vec{v_j} = c^{(i)}_1 \vec{v_1} + \cdots c^{(i)}_n\vec{v_n}, \hspace{1cm} c^{(i)}_j \in \mathbb R,~ i =1,\dots m.
\end{align} 

\noindent where each $\vec v_j$ is the $j^{th}$ singular vector of $A$. In the context of our problem, it should not be assumed that knowledge of the singular vectors be known; otherwise, the solution could be reconstructed and the problem would be trivial. However, we analyze the convergence using this knowledge to inspire a new algorithm. Considering that $A = U\Sigma V^{\mathsf T}$, we obtain each $c^{(i)}_n$ by:

\begin{align}
    c_n= [c^{(1)}_n, \dots, c_n^{(m)}]^{\mathsf T} = AV_{:n}
\end{align}
With these $n$th spectral coefficients $c_n$, we construct spectral sampling weights $\omega_i$ for each row of the system:
\begin{align}\label{spectral_coeffs}
    \omega_i := \frac{|c^{(i)}_n|}{\sum_{l=1}^m|c^{(l)}_n|}, \hspace{1cm} i=1,\dots,m.
\end{align}

\noindent By substituting the row norm-based sampling distribution by the one determined by the spectral weights from (\ref{spectral_coeffs}), we run a series of numerical experiments, expecting to improve upon the approximation error horizon encountered by the RK algorithm.

\subsection{Numerical Experimentation}\label{sub:numericalExperimentation}

For all numerical experiments in this section, we refer to the following matrix construction: We produce $A \in \mathbb R^{240 \times 12}$ by taking the QR-decompositions of the Gaussian matrices $\hat U \in \mathbb R^{240 \times 240}$ and $\hat V \in \mathbb R^{12 \times 12}$ to obtain orthogonal matrices $U \in \mathbb R^{240 \times 240}$ and $V \in \mathbb R^{12 \times 12}$. We then construct a matrix $\Sigma \in \mathbb R^{240 \times 12}$ with entries $\Sigma_{jj}=e^{n-j+1} , ~(j = 1,\dots, n)$ and 0 everywhere else such that we obtain $A = U\Sigma V^T$. To complete the system, we let $x^* = [0,0,\dots, 1]^{\mathsf T}$ be the length-$n$ solution vector to $Ax^* = b$.\\

\noindent To quantify the progress of the iterate $x_k$ in the direction of each singular vector, we follow the dot product of the error vector $\varepsilon_k$ with each singular vector. We define the $j$th singular error at each iteration as:

\begin{align}\label{singularErr}
    \varepsilon_k^{(j)} := |\langle x-x^*, v_j\rangle|, \hspace{1cm}j=1,\dots, n
\end{align}

\begin{figure}[h!]
    \centering
    \includegraphics[width=1\linewidth]{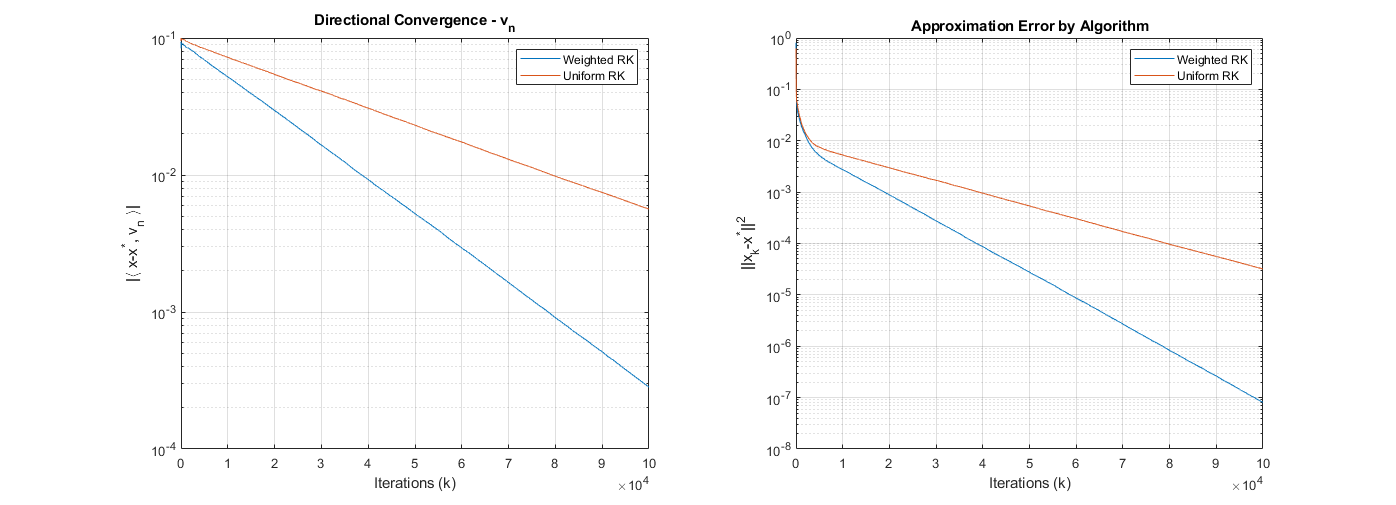}
    \caption{A matrix $A \in \mathbb R^{240 \times 12}$ is constructed as described in Section \ref{sub:numericalExperimentation}. Plot data averages results over 50 different random matrices. The left plot describes the directional convergence $|\langle x_k-x^*, v_n\rangle|$ for each algorithm. The right plot demonstrates the improved convergence rate for weighted sampling. }
    \label{fig:directionalConvergence}
\end{figure}

\noindent We notice in Figure \ref{fig:directionalConvergence} that the modified sampling method inspired by the directional weighting scheme improves the convergence in the direction of the least singular vector. It does not improve the directional convergence in any of the other singular vectors. However, this is expected as we determine our sampling weights using only the coefficient $c_n$. This method effectively targets the restricting factor of (\ref{expectedConvergence}) by  sampling each vector with a probability proportional to the  magnitude of its component in the direction of $v_j$. Figure \ref{fig:directionalConvergence} also reveals that this weighted sampling distribution, on average, improves the overall convergence rate of the algorithm. Further plotting would reveal that $x_k$ converges to $x^*$ within machine error $(\epsilon_{machine} = 1\times10^{-15})$ in less than half the number of iterations when using weighted sampling compared to uniform.\\

\noindent Although this method requires knowledge of the singular value decomposition of $A$, the trend it unveils should be used to inspire methods that improve existing RK algorithms for ill-conditioned systems. We propose some solutions that leverage these properties in Section \ref{futureWork}.

\section{Future Work}\label{futureWork}

We remind the reader that the weighted sampling method proposed in Section \ref{sec:singular_directions} requires the knowledge of the SVD of the matrix beforehand. Since the process of obtaining the SVD is, itself, expensive and would render the problem trivial, we suggest approximating the singular vectors $v_j$ of $A$ iteratively. We suggest two possible approaches: (i) iteratively adapt sampling weights as the singular vectors are approximated, or (ii) perform Average Block Kaczmarz (ABK) \cite{Kui2020} with $\beta$ many rows and adjust the projection weight of each row proportionally to the magnitude of its component $c^{(i)}_n$ relative to the other rows'.

\newpage

\end{document}